\documentclass[11pt]{article}
\usepackage{amsfonts,color}
\usepackage{amsmath,amssymb,tipa,enumerate,txfonts,mathrsfs,algorithm,float,indentfirst}
\usepackage{bm}
\usepackage{a4wide}
\usepackage{epsfig,multirow,epstopdf}
\usepackage{graphicx}

\makeatletter
  \newcommand\figcaption{\def\@captype{figure}\caption}
  \newcommand\tabcaption{\def\@captype{table}\caption}
\makeatletter

\newtheorem{lemma}{Lemma}[section]
\newtheorem{theorem}{Theorem}[section]

\def\be{\begin{equation}}
\def\en{\end{equation}}
\def\beq{\begin{eqnarray}}
\def\eq{\end{eqnarray}}
\def\beqx{\begin{eqnarray*}}
\def\eqx{\end{eqnarray*}}

\def\12{{1\over 2}}

\def\0{{\bf 0}}

\def\prod{{\bf \sqcap}}

\def\Pi{\pi}

\begin{document}
\title{Discontinuous hybrid neural networks for the one-dimensional partial differential equations}
\author{
Xiaoyu Wang\footnote{College of Mathematics and Systems Science, Shandong
University of Science and Technology,
Qingdao 266590, China (2836401677@qq.com).} 
~~~Long Yuan\footnote{Corresponding author. College of Mathematics and Systems Science, Shandong University of Science and Technology,
Qingdao 266590, China. The author was supported by Shandong Provincial Natural Science Foundation under the grant ZR2024MA059. (email:\it{sdbjbjsd@163.com}). 
The author proposed the ideas of the method and analysis. 
}
~~~and~~~
Yao Yu\footnote{
College of Mathematics and Systems Science, Shandong University of Science and Technology, Qingdao 266590, China. The author was supported by Shandong Provincial Natural Science Foundation under the grant ZR2021QA109. (yuyao0608@163.com).}
}

\date{\today }
\maketitle

\begin{abstract} A feedforward neural network, including hidden layers, motivated by nonlinear functions (such as Tanh, ReLU, and Sigmoid functions), exhibits uniform approximation properties in Sobolev space, and discontinuous neural networks can reduce computational complexity. In this work, we present a discontinuous hybrid neural network method for solving the partial differential equations, construct a new hybrid loss functional that incorporates the variational of the approximation equation, interface jump stencil and boundary constraints. The RMSprop algorithm and discontinuous Galerkin method are employed to update the nonlinear parameters and linear parameters in neural networks, respectively. This approach guarantees the convergence of the loss functional and provides an approximate solution with high accuracy.
\end{abstract}


\vskip 0.1in


\section{Introduction}

In numerical methods for solving partial differential equations, neural networks possess excellent adaptive capabilities and have gradually become a powerful tool\cite{maochao}. A Physics-Informed Neural Network(PINN)\cite{ref1}, as a classical neural network method, embeds the physical laws or boundary conditions directly into the training process of the neural network, which ensures that the training results conform to the physical laws. This approach can be used for control optimization, and forward and inverse problems \cite{Cai}. Compared to purely data-driven neural networks, PINNs learn more generalized models with fewer data samples. Furthermore, \cite{Zhang} proposed a Variational Physics-Informed Neural Network(VPINN) within the Petrov-Galerkin framework. Unlike the PINN that incorporates the strong form of the equation into the network, a VPINN incorporates the variational weak formulation of the problem and constructs a variational formulation of loss function. While such PINNs and VPINN approaches have enjoyed success in particular cases, in practice, the achievable relative error is prone to stagnate around 0.1-0.01 no matter how many neurons or layers are uesd to define the underlying network architecture (see \cite{Jagtap}).

To overcome this plight, \cite{Ainsworth} proposed Galerkin neural networks to approximate the variational equations by adaptively constructing a sequence of finite dimensional neural network subspaces, which are used as the discrete space of the variational equations. It constructs discrete basis functions under Galerkin iteration and provides a novel framework for solving elliptic, parabolic, and hyperbolic partial differential equations. However, the approach utilizes continuous neural networks, which typically require a large number of test points to ensure accuracy, thereby increasing computational complexity.

In this paper, we use discontinuous neural networks to solve boundary value problems for time-harmonic partial differential equations. For convenience, we call the proposed method as DHNN method. For each element, fewer neural network nodes are used to construct the network structure, which better adapts to the local features and reduces the requirement of test point density\cite{Yuanhu}. Inspired by the construction of loss functionals in variational physics-informed neural networks\cite{Zhang}, we propose an innovative hybrid loss functional in section 3 that incorporates the variational of the approximation equation, interface jump stencil and boundary constraints. The RMSprop optimizer is employed to search nonlinear network parameters, such as weights and biases. After updating the nonlinear parameters, based on the associated linear neural network space, the discontinuous Galerkin method is employed to solve the linear parameters in the neural network. We derive the weak convergence of the approximation solution in the sense of the loss functional in section 4. Finally, we report some numerical results to confirm the effectiveness of the proposed method in section 5.

\section{The model and discontinuous neural networks}
\subsection{The model}
We consider the solution of boundary value problems for one-dimensional time-harmonic partial differential equations using discontinuous neural networks. However, the proposed method can be extended to two-dimensional and three-dimensional cases. Let $\Omega\subset R$ be a bounded domain, we perform a mesh discretization of the solution domain, with nodes being:
\begin{equation*}
x_0<x_1<\cdots<x_N.
\end{equation*}
Here the subdomain $\Omega_i=[x_{i-1},x_i]$ between the adjacent nodes $x_{i-1}$ and $x_i$ is $i$-th element. Let ${\cal T}_h$ denote the partitioning $\{\Omega_k\}_{k=1}^{N}$ of domain $\Omega$.

Consider the following linear problem
\begin{equation} \label{variahelm1}
\left\{ \begin{aligned}
     & \mathcal{A}u(x) = f  \quad \quad\quad\text{in} \quad \Omega ,\\
    & \mathcal{B}u(x) = g(x)  \quad\quad \text{on} \quad\partial\Omega ,
                          \end{aligned} \right.
                          \end{equation}
$f\in L^2(\Omega)$, $g\in H^{\frac{1}{2}}(\partial\Omega)$. Let $(V({\cal T}_h), ||\cdot||_V)$, $(Y, ||\cdot||_Y)$, $(Z, ||\cdot||_Z)$ and $(W, ||\cdot||_W)$ be Banach spaces. Suppose that $\mathcal{A}:V({\cal T}_h)\longrightarrow Y$ is a bounded linear operator and $\mathcal{B}:V({\cal T}_h)\longrightarrow Z$ is also a bounded linear operator. Besides, suppose that $C_i(i=1,2):V({\cal T}_h)\longrightarrow W$ are bounded linear operators related to jump quantities  on the set of all interior points $\{x_i\}_{i=1}^{N-1}$.


\subsection{Discontinuous neural networks}
A general discontinuous neural network consists of a single hidden layer of $n\in\mathbb{N}$ neurons on each elemen $\Omega_k\in{\cal T}_h$, defining a function $U_{NN}:\mathbb{R}\longrightarrow\mathbb{C}$ as follows:
\begin{equation}
u_{NN}(x;\theta)|_{\Omega_k}=\sum_{j=1}^{n}c_j^{(k)}\sigma(W_j^{(k)}\cdot x+b_j^{(k)}) ,\quad\quad\quad\forall\Omega_k\in V({\cal T}_h),
\end{equation}
where $n$ is called as the width of the discontinuous network; $W_j^{(k)}\in\mathbb{R}$, $b_j^{(k)}\in\mathbb{R}$ are elementwise nonlinear parameters; $c_j^{(k)}\in\mathbb{C}$ are elementwise linear parameters; and $\sigma:\mathbb{R}\longrightarrow\mathbb{C}$ is a bounded, elementwise smooth activation function. For simplicity, we shall use the elementwise notation $W|_{\Omega_k}=[W_1^{(k)},\cdots,W_n^{(k)}]\in\mathbb{R}^{1\times n}$ referred as weights, $b|_{\Omega_k}=[b_1^{(k)},\cdots,b_n^{(k)}]\in\mathbb{R}^n$ referred as biases, and $c|_{\Omega_k}=[c_1^{(k)},\cdots,c_n^{(k)}]^T\in\mathbb{C}^n$. The set of nonlinear and linear parameters is collectively denoted by $\Theta=[W,b,c]$.  In particular, a set of nonlinear parameters defined on every element is collectively denoted by $\Phi=\{W,b\}$.

Denote $V_{n}^{\sigma}({\cal T}_h)$  as the set of all functions of the form(2.2), which is defined by
\begin{equation*}
V_{n}^{\sigma}({\cal T}_h):=\{v:v|_{\Omega_k}=\sum_{j=1}^{n}c_j^{(k)}\sigma(W_j^{(k)}\cdot x+b_j^{(k)}),\quad k=1,\cdots, N\}.
\end{equation*}


Although neural network space $V_{n}^{\sigma}({\cal T}_h)$ is neither closed nor compact in $C(\Omega)$ or $L^p(\Omega)$, by restricting the boundedness of network parameters, the neural network space with bounded parameters is both closed and compact in $C(\Omega)$ and $L^p(\Omega)$.  For the well-posedness of the discretized optimization problem (3.3) defined below, the subset $V_{n,C}^{\sigma}({\cal T}_h)$ of $V_{n}^{\sigma}({\cal T}_h)$ consisting of functions with bounded parameters needs to be defined (see \cite{Ainsworth}).

Given a sequence of positive increasing numbers $C_n$, satisfying$ \underset{n\longrightarrow\infty}{\lim}C_n=\infty$, we define
\begin{equation*}
V_{n,C}^{\sigma}({\cal T}_h):=\{v\in V_{n}^{\sigma}({\cal T}_h):||\Theta||<C_n\},
\end{equation*}
where$||\Theta||=max_k||\Theta_k||, ||\Theta_k||:= \max\limits_j|W_{j}^{(k)}|+ \max\limits_j|b_j^{(k)}|+ \max\limits_j|c_j^{(k)}|$. 

\section{A discontinuous hybrid neural network algorithm}
Suppose $\{\varphi_{k,i},k=1,\cdots,N,i=1,\cdots,\infty\}$ forms a complete orthonormal basis in a Hilbert space Y. Also, let $\{\varphi_{k,i}\}_{i=1}^{\infty}$ forms a complete set of orthogonal basis functions in the Hilbert space $Y_k=Y|_{\Omega_k}$ . We denote the PDE residual loss by $L_{r}^{h}(\Theta):=\sum\limits_{k=1}^{N}\sum\limits_{i=1}^{\infty}|(\mathcal{A}v-f,\varphi_{k,i})_{Y_k}|^{2}$ , the boundary residual loss by $L_b(\Theta):=||\mathcal{B}v-g||_{Z}^{2}$ , and the interface residual loss $L_{int}^{h}(\Theta):=||[[C_1v]]||_{W}^{2}+||[[C_2v]]||_{W}^{2}$ . Then we can define the continuous hybrid loss functional as follows:
\begin{equation}
\mathcal{J}_{\tau}^{h}(v)=L_{r}^{h}(\Theta)+\tau(L_b(\Theta)+L_{int}^{h}(\Theta)),     \quad\quad\quad\forall\tau>0,v\in V({\cal T}_h).
\end{equation}
Correspondingly, define its corresponding truncation by $L_{r}^{h,\mathbf{M}}(\Theta):=\sum\limits_{k=1}^{N}\sum\limits_{i=1}^{M_k}|(\mathcal{A}v-f,\varphi_{k,i})_{Y_k}|^{2}$ and further
\begin{equation}
\mathcal{J}_{\tau}^{h,\mathbf{M}}(v)=L_{r}^{h,\mathbf{M}}(\Theta)+\tau(L_b(\Theta)+L_{int}^{h}(\Theta)), \quad\quad\quad\forall\tau>0,v\in V({\cal T}_h),
\end{equation}
where $\mathbf{M}=(M_1,\cdots,M_N)$.

Now we consider the minimization problem of the loss functionals constrained on the discontinuous neural network space. That is
\begin{equation} \label{minpro}
\underset{v\in V_{n,C_n}^{\sigma}}{\min}\mathcal{J}_{\tau}^{h,\mathbf{M}}(v).
\end{equation}

To calculate optima of the objective function $\mathcal{J}_{\tau}^{h,\mathbf{M}}(v)$ in neural network trainings, the RMSprop algorithm and discontinuous Galerkin method are employed to update the nonlinear parameters and linear parameters, respectively. The learning rate annealing strategy \cite{Wangteng} is used to adaptively update the $\tau$ for  balancing the the PDE residual loss $L_{r}^{h,\mathbf{M}}(\Theta)$ and the boundary residual loss $L_b(\Theta)+L_{int}^{h}(\Theta)$.

\subsection{Updating linear parameters}
For known nonlinear parameters $\Phi^l$ with $l\geq0$, let $V_{\Phi^l}({\cal T}_h)$ denote the space $V_{n}^{\sigma}({\cal T}_h)$ with the given nonlinear parameters $\Phi^l$. Linear parameters $c^l$ are obtained using the standard discontinuous Galerkin method. Define the sesquilinear form $a(\cdot,\cdot)$ by
\beq
a(u,v) &= &\sum_{k=1}^{N}\sum_{i=1}^{M_k}(\mathcal{A}u,\varphi_{k,i})_{Y_k}\overline{(\mathcal{A}v, \varphi_{k,i})_{Y_k}}
+ \tau \big( (\mathcal{B}u,\mathcal{B}v)_z
\cr & + &
([[C_1u]],[[C_1v]])_Z + ([[C_2u]],[[C_2v]])_Z \big),\quad\quad\quad\quad\forall u,v\in V_{\Phi^l}({\cal T}_h),
\eq
and semilinear form $L(\cdot)$ by
\begin{equation}
L(v)=\tau(g,\mathcal{B}v)_z+\sum_{k=1}^{N}\sum_{i=1}^{M_k}(f,\varphi_{k,i})_{Y_k}\overline{(\mathcal{A}v,
\varphi_{k,i})_{Y_k}}, \quad\quad\quad\quad\forall v\in V_{\Phi^l}({\cal T}_h).
\end{equation}
Then  the minimization problem (\ref{minpro}) can be expressed as follows: find $u_l\in V_{\Phi^l}({\cal T}_h)$ such that
\begin{equation}
a(u_l,v) =L(v), \quad\quad\quad\forall v\in V_{\Phi^l}({\cal T}_h).
\end{equation}
\subsection{Updating nonlinear parameters}
For known nonlinear parameters $c^l$ with $l>0, V_{c^l}({\cal T}_h)$ denotes the set of neural network functions under the corresponding activation coefficients. The nonlinear parameters $\Phi^{l+1}$ are updated using the stochastic gradient descent algorithm RMSprop\cite{Mada}, satisfying
\begin{equation}
\mathcal{J}_{r}^{h,\mathbf{M}}(\tilde{u}_{l+1})=\underset{v\in V_{c^l}({\cal T}_h)}{\min}\mathcal{J}_{\tau}^{h,\mathbf{M}}(v).
\end{equation}

\subsection{Updating the loss balancing factor $\tau$}
In order to balance the interplay between the different loss terms in the loss function, the following formula is used to update the balancing factor $\tau$ (see \cite{Wangteng}). Define
\begin{equation}
\hat{\tau}=\frac{\max_{\Theta_l}{\{|\bigtriangledown_\Theta L_{r}^{h,\mathbf{M}}(\Theta_l)|\}}}{\overline{\overline{|\bigtriangledown_\Theta(L_b(\Theta_l)+L_{int}^{h}(\Theta_l))|}}}.
\end{equation}
where $\Theta_l$ denotes the values of the parameters $\Theta$ at the $l-$th iteration. $|\cdot|$ denotes the elementwise absolute value, and $\overline{\overline{|\bigtriangledown_\Theta(L_b(\Theta_l)+L_{int}^{h}(\Theta_l))|}}$ denotes the mean of $|\bigtriangledown_\Theta(L_b(\cdot)+L_{int}^{h}(\cdot))|$ at $\Theta_l$. Thus, in the $(l+1)$-th iteration, the adaptive iteration formula for $\tau$ is given by
\begin{equation}
\tau=(1-\beta)\tau+\beta\hat{\tau}.
\end{equation}
The recommended hyperparameter $\beta$ is set to be $0.1$.

\subsection{Algorithm overview}
This section presents the implementation of the Discontinuous Neural Network algorithm:\\

\textbf{Step 1} initialize input parameters: domain $\Omega$, $\omega$, tolerance $\rho$, nonlinear parameters $\Phi_0$;

\textbf{Step 2} Update the linear parameters using (3.6);

\textbf{Step 3} Update $\tau$ using (3.8)-(3.9);

\textbf{Step 4} Use the RMSProp algorithm\cite{Wangteng} to solve (3.7) for updating the nonlinear parameter $\Phi^{l+1}$;

\textbf{Step 5} If $||\Phi_{(l+1)}-\Phi_l||<\rho$, the iteration is terminated, else go to \textbf{Step 2}. \\

In particular,  we call the process of updating the nonlinear network parameters using \textbf{Step 4} as inner iteration. Here, the RMSprop optimizer utilizes an initial learning rate, denoted as $\alpha$, which determines the starting step size for updating model parameters during the neural network trainings. The initial learning rate may vary across different models to accommodate their unique characteristics and requirements.

\section{Convergence estimates}
With the help of classical approximation process\cite{Hornik} , discontinuous neural networks have the following the universal approximation property.
\begin{lemma}
Suppose that $1\leq p <\infty$, $0\leq s <\infty$, and that $\Omega\subset R^d$ is compact, then $V_{n}^{\sigma}({\cal T}_h)$ is dense in $W^{s,p}({\cal T}_h):=\{v\in L^p(\Omega), D^\alpha v\in L^p(\Omega_k) \quad \forall |\alpha|\leq s,k=1,\cdots N\}$. In particular, for any given function $f\in W^{s,p}({\cal T}_h)$ and $\tau>0$, there exists $n(\tau,f)\in\mathbb{N}$ and $\tilde{f}\in V_{n(\tau,f)}({\cal T}_h)$ such that $||f-\tilde{f}||<\tau$.
\end{lemma}

Let $\{\varphi_{k,i}\}_{i=1}^{M_k}$ be a set of orthogonal basis functions in the Hilbert space $Y_k$. With this basis, define the orthogonal projection operator as $P_{h,\mathbf{M}}$: $Y\longrightarrow Y_k$
\begin{equation}
(P_{h,\mathbf{M}}v)|_{\Omega_k}:=\sum_{i=1}^{M_k}(v,\varphi_{k,i})_{Y_k}\varphi_{k,i}.
\end{equation}

By the orthogonality, we have
\begin{equation}
\mathcal{J}_{\tau}^{h,\mathbf{M}}(v)=\mathcal{J}_{\tau}^{h}(v)-\sum_{k=1}^{N}||(I-P_{h,\mathbf{M}})(\mathcal{A}v-f)
||_{Y_k}^{2},\quad\quad\quad\quad\forall v\in V(\mathcal{T}_h).
\end{equation}
By the universal approximation property, we obtain the convergence of discontinuous hybrid neural network approximation in the sense of loss functional.

\begin{theorem}
For a fixed $\tau>0$ and $\mathbf{M}$, it holds that $\mathcal{J}_{\tau}^{h,\mathbf{M}}(u_{n}^{\tau,\mathbf{M}})\longrightarrow0$ when $n\longrightarrow\infty$.
\end{theorem}

\textit{Proof}. By universal approximation property, we have $\underset{v\in V({\cal T}_h)}{\min}\mathcal{J}_{\tau}^{h}(v)=0$. By (4.2), we have $\underset{v\in V({\cal T}_h)}{\min}\mathcal{J}_{\tau}^{h,\mathbf{M}}(v)=0$. Since $\underset{n}{\cup}V_{n,C_n}^\sigma=V({\cal T}_h)$ and $V_{n,C_n}^\sigma\subset V_{{n+1},C_{n+1}}^\sigma$, by letting $n\longrightarrow\infty$, the proof is completed. $\hfill\Box$

Additionally, we can obtain the following an a posterior error estimate.
\begin{theorem}
For any $\epsilon>0$, there exists $\mathbf{M}_\epsilon$ such that, for all $v\in V_{n,C_n}^\sigma$, we have
\begin{equation*}
\mathcal{J}_{\tau}^{h}(v)-\epsilon\leq\mathcal{J}_{\tau}^{h,\mathbf{M}_\epsilon}(v)\leq\mathcal{J}_{\tau}^{h}(v),
\end{equation*}
where $\mathbf{M}_\epsilon=(\mathbf{M}_{\epsilon,1},\cdots,\mathbf{M}_{\epsilon,N})$.
\end{theorem}

\textit{Proof}. Considering that $V_{n,C_n}^\sigma$ is a compact set in $V({\cal T}_h)$, and $\mathcal{A}:V({\cal T}_h)\longrightarrow Y$ is a bounded linear operator , so the set $\mathcal{A}V_{n,C_n}^\sigma$ is still compact in $Y$. Thus, for any $\epsilon>0$, there exists a finite dimensional subspace $\tilde{Y}_{\mathbf{M}_\epsilon}\subset Y$ such that for any $w\in\mathcal{A}V_{n,C_n}^\sigma-f$, there exists $w_\epsilon\in\tilde{Y}_{\mathbf{M}_\epsilon}$ satisfying $||w_\epsilon-w||_Y\leq\epsilon^\frac{1}{2}$. Set $\tilde{Y}_{\mathbf{M}_\epsilon}=span\{\varphi_{k,i}:i=1,\cdots,\mathbf{M}_{\epsilon,k}, k=1,\cdots,N\}$. Combining with (4.2), we have
\be
\mathcal{J}_{\tau}^{h}(v)-\epsilon\leq\mathcal{J}_{\tau}^{h,\mathbf{M}_\epsilon}(v)\leq\mathcal{J}_{\tau}^{h}(v). 
\en $\hfill\Box$

\section{Numerical examples}
In this section we report some numerical results to confirm the effectiveness of the proposed method. For comparison purposes, we also implemented the VPINN approach, employing the Xavier initialization scheme for its nonlinear parameters. To calculate optima of the objective function $\mathcal{J}_{\tau}^{h,\mathbf{M}}(v)$, we employ a fixed Gauss-Legendre quadrature rule with 30 nodes to approximate all inner product. Let $\Omega=[0,1]$. In the discontinuous hybrid neural network, two different activation functions are used: the sigmoid function: $\sigma(x)=\frac{1}{1+e^{-x}}$, and the hyperbolic tangent function: $tanh(x)=\frac{e^x+e^{-x}}{e^x-e^{-x}}$. We set the tolerance $\rho=10^{-6}$ and  the inner iteration count to be 50, meaning that after every 50 internal iterations for updating the nonlinear parameters $W$ and $b$ 50, there is one iteration for updating the linear parameter $c$. All calculations are implemented on python.

Define jump operators on the set of all interior points:
\begin{gather*}
[[C_1u]]=\omega^2[[u]]=\omega^2(u_k-u_j),    \\
[[C_2u]]=[[\partial_{n}u]]=\partial_{nk}u_k+\partial_{nj}u_j.
\end{gather*}

\subsection{Possion equations}
Consider Possion equations  which is formalized,  by
\begin{equation}
\left\{ \begin{aligned}
     & -\Delta u(x) = f  \quad \quad\quad \quad\quad \quad\quad\text{in} \quad \Omega ,\\
    &   \partial_n u(x)+u(x) = g   \quad\quad\quad\quad\quad\text{on} \quad\partial\Omega,
                          \end{aligned} \right.
                          \end{equation}
where $g\in H^{\frac{1}{2}}(\partial\Omega)$, $f\in L^2(\Omega)$. The outer normal derivative is referred to by $\partial_n$ . In the numerical solution of the Poisson equation, the relaxation parameter $\omega$ in the jump operator is taken as $2\pi$.
\subsubsection{ A homogeneous case}
Consider homogeneous Possion equation, i.e. $f=0$. Let the exact solution be $u=2x+1$. The initial values of all nonlinear parameters are set to be 0 in DHNN.

\begin{figure}[H]
\begin{center}
\begin{tabular}{cc}
 \epsfxsize=0.5\textwidth\epsffile{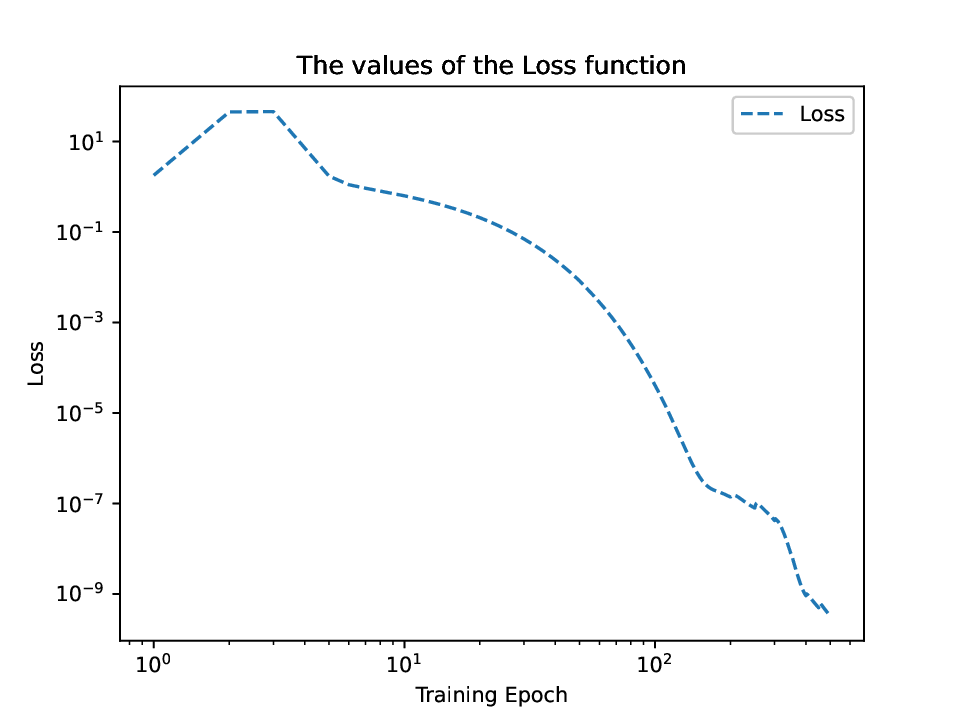}&
 \epsfxsize=0.5\textwidth\epsffile{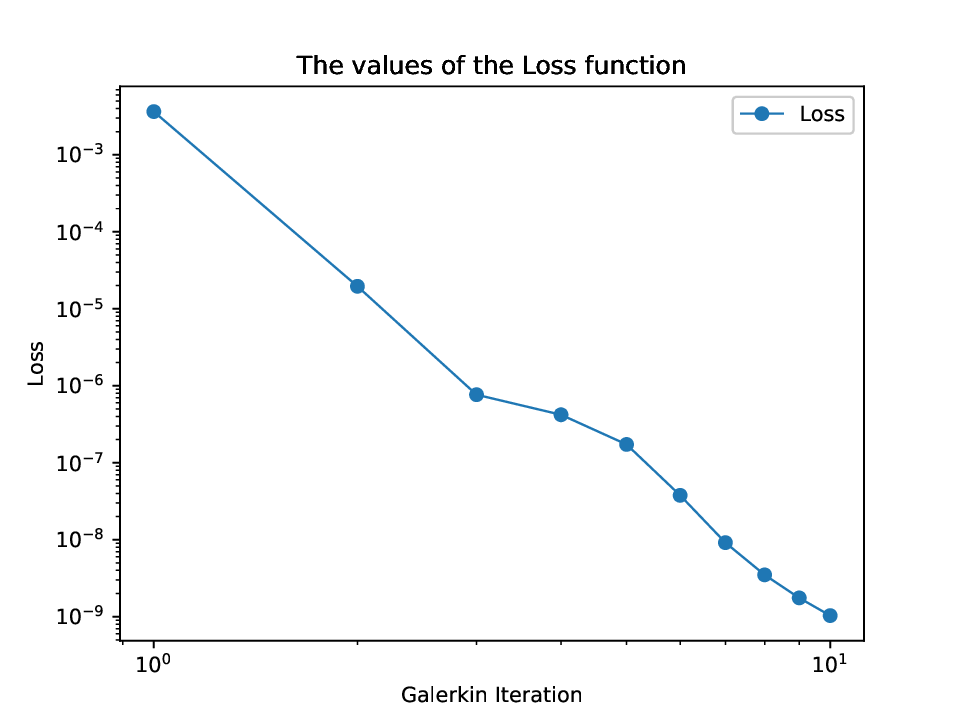}\\
\end{tabular}
\end{center}
 \caption{ The values of the loss functional for the DHNN. Left: inner iterations. Right: outer iterations.}
\label{Possion_loss_sig}
\end{figure}

Figure \ref{Possion_loss_sig} shows the optimization process of the loss functional for both inner and outer iterations, with the activation function of the hidden layer being the sigmoid function. The parameters are set  as follows:
\begin{eqnarray*}
h=\frac{1}{N}=\frac{1}{4},n=10, \alpha=0.04.
\end{eqnarray*}

Figure \ref{Possion_loss_sig} exhibits excellent convergence behavior of the loss functional, achieving a six-order-of-magnitude reduction within just 200 iterations. Such rapid convergence demonstrates the remarkable efficiency of the proposed optimization algorithm.

\begin{figure}[H]
\begin{center}
\begin{tabular}{cc}
 \epsfxsize=0.5\textwidth\epsffile{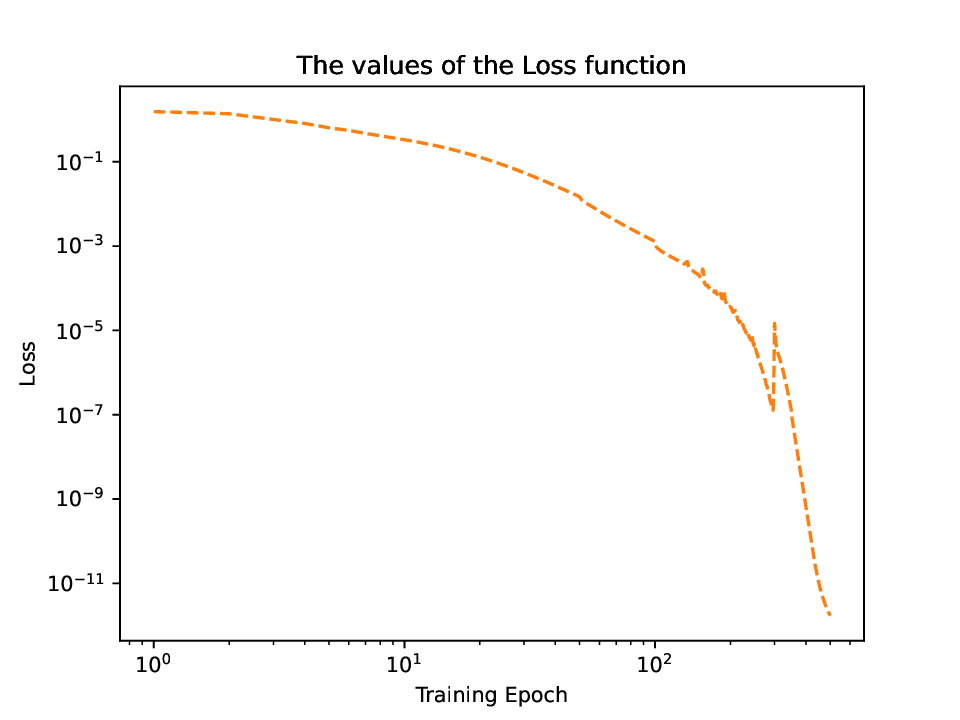}&
 \epsfxsize=0.5\textwidth\epsffile{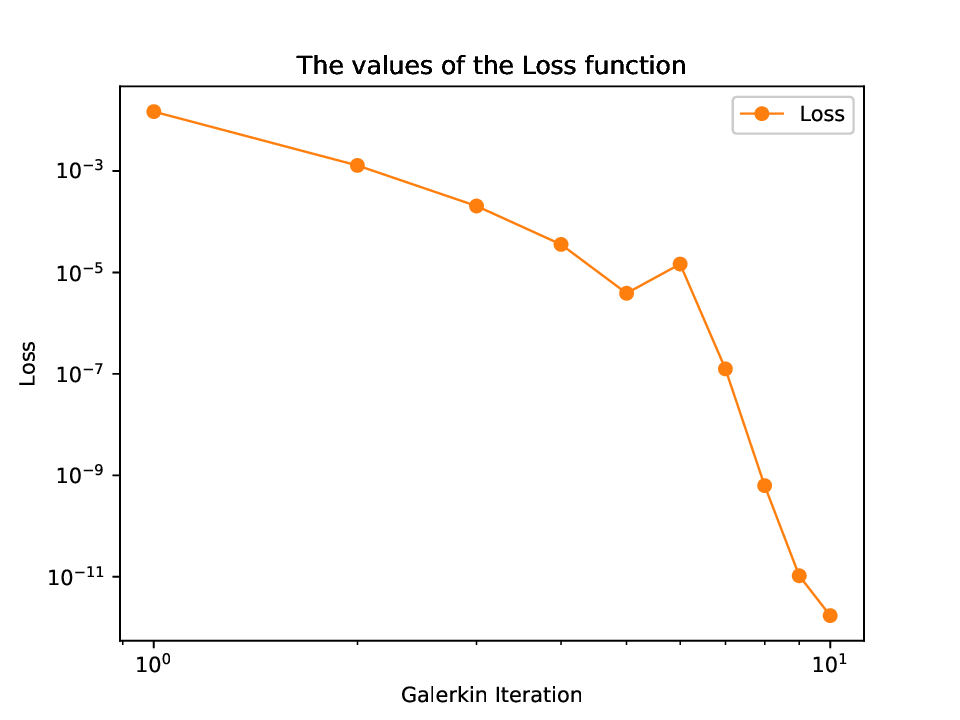}\\
\end{tabular}
\end{center}
 \caption{ The values of the loss functional for the DHNN. Left: inner iterations. Right: outer iterations.}
\label{Possion_loss_tanh}
\end{figure}

Figure \ref{Possion_loss_tanh} shows the optimization process of the loss functional for both inner and outer iterations, with the activation function of the hidden layer being the tanh function. The parameters are set as follows:
\begin{eqnarray*}
h=\frac{1}{N}=\frac{1}{4},n=10, \alpha=0.05.
\end{eqnarray*}
As shown in Figure \ref{Possion_loss_tanh}, the loss function demonstrates an overall favorable descending trend, ultimately achieving an eleven-order-of-magnitude reduction. However, distinct oscillations are observed during the optimization process, where the loss value abruptly increases by two orders of magnitude within merely two iterations. Notably, the post-oscillation convergence rate becomes significantly faster than the pre-oscillation phase.

\begin{figure}[H]
\begin{center}
\begin{tabular}{cc}
 \epsfxsize=0.5\textwidth\epsffile{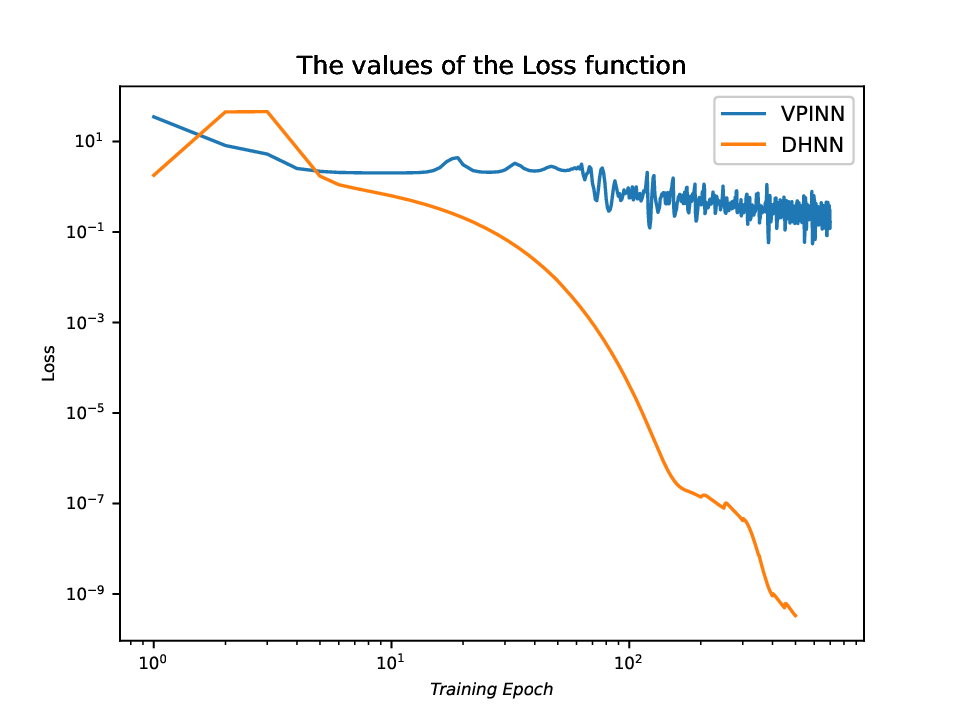}&
 \epsfxsize=0.5\textwidth\epsffile{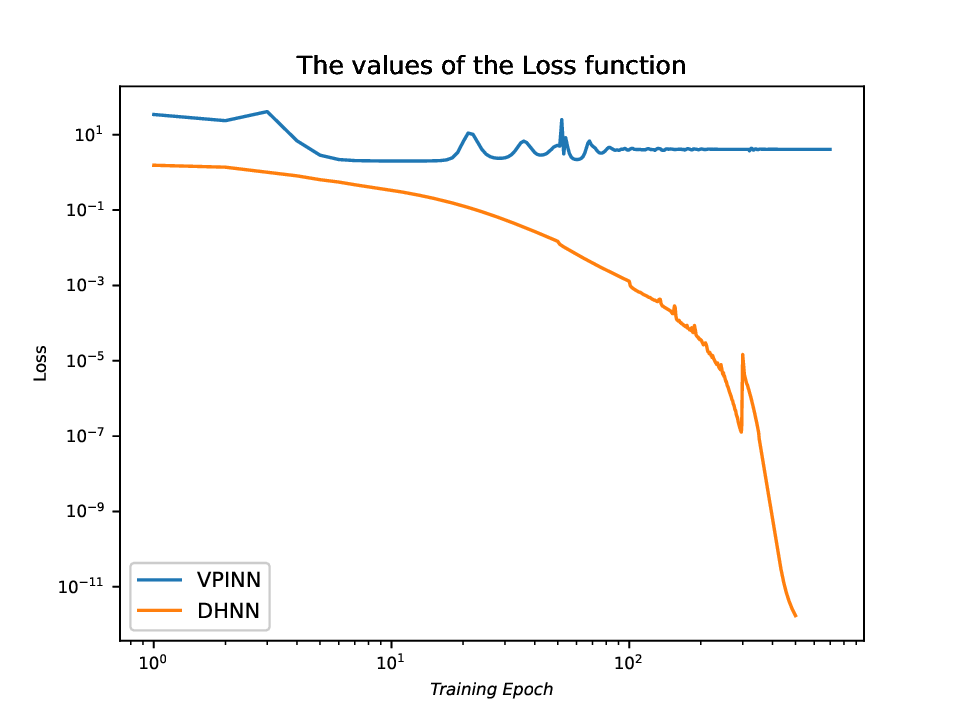}\\
\end{tabular}
\end{center}
 \caption{  Comparison of the loss functional values between the VPINN and the DHNN methods. Left: sigmoid activation function. Right: tanh activation function.}
\label{Possion_loss_VPINN}
\end{figure}

Figure \ref{Possion_loss_VPINN} shows the loss optimization for both VPINN and DHNN methods with sigmoid and tanh activation functions. The parameters of the VPINN method are set as follows:
\begin{align*}
\text{sigmoid}: &\quad\text{four hidden layers (40 neurons per layer)},\ \alpha=0.05;\\
\text{tanh}:    &\quad\text{four hidden layers (40 neurons per layer)},\ \alpha=0.05.
\end{align*}
Figure \ref{Possion_loss_VPINN} indicates that, regardless of the activation function employed, the VPINN method fails to achieve effective convergence after 500 iterations. The optimization process under both activation functions exhibits significant oscillatory behavior without observable amplitude attenuation. In contrast, the DHNN method demonstrates superior convergence characteristics under identical conditions.

\begin{figure}[H]
\begin{center}
\begin{tabular}{cc}
 \epsfxsize=0.5\textwidth\epsffile{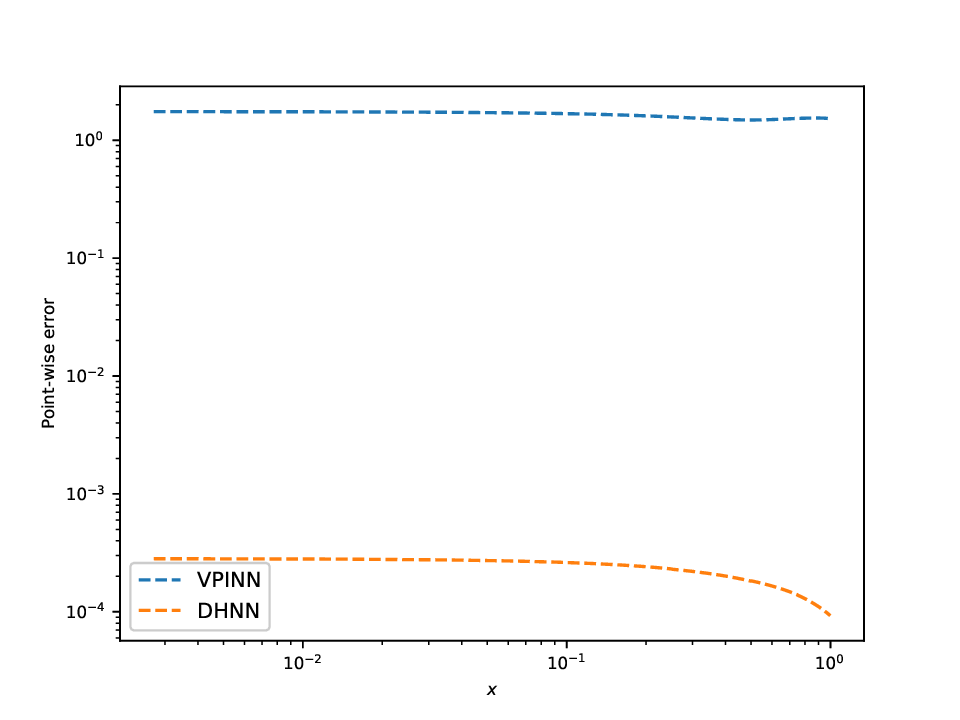}&
 \epsfxsize=0.5\textwidth\epsffile{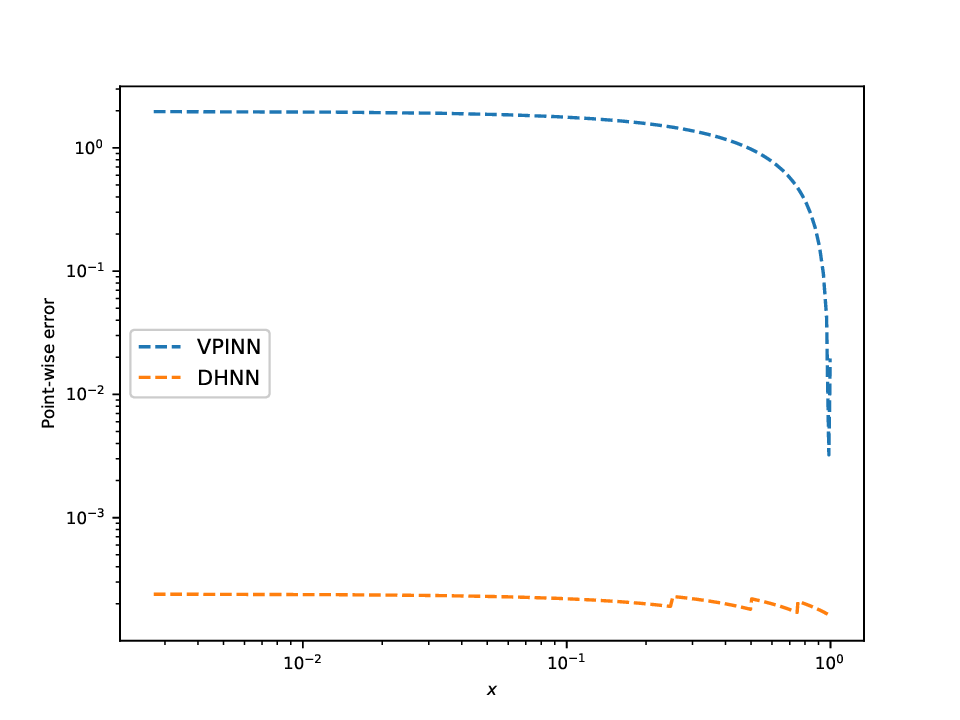}\\
\end{tabular}
\end{center}
 \caption{ The point-wise error distributions for both methods. Left: sigmoid activation function. Right: tanh activation function.}
\label{Possion_error_compare}
\end{figure}

Figure \ref{Possion_error_compare} presents the point-wise error distributions for solving the homogeneous Poisson equation using both the VPINN method and the DHNN method with sigmoid and tanh activation functions.
Figure \ref{Possion_error_compare} demonstrates that the VPINN method fails to achieve effective error reduction or satisfactory numerical accuracy across different activation functions. In contrast, the DHNN algorithm exhibits three distinct advantages: consistently maintaining errors at $10^{-4}$ magnitude, displaying more uniform error distribution with faster convergence rates, and preserving stable high-precision characteristics regardless of activation function. These results confirm DHNN's superior adaptability, outperforming VPINN in all tested configurations.
\subsubsection{An inhomogeneous case}
\begin{figure}[H]
\begin{center}
\begin{tabular}{cc}
 \epsfxsize=0.5\textwidth\epsffile{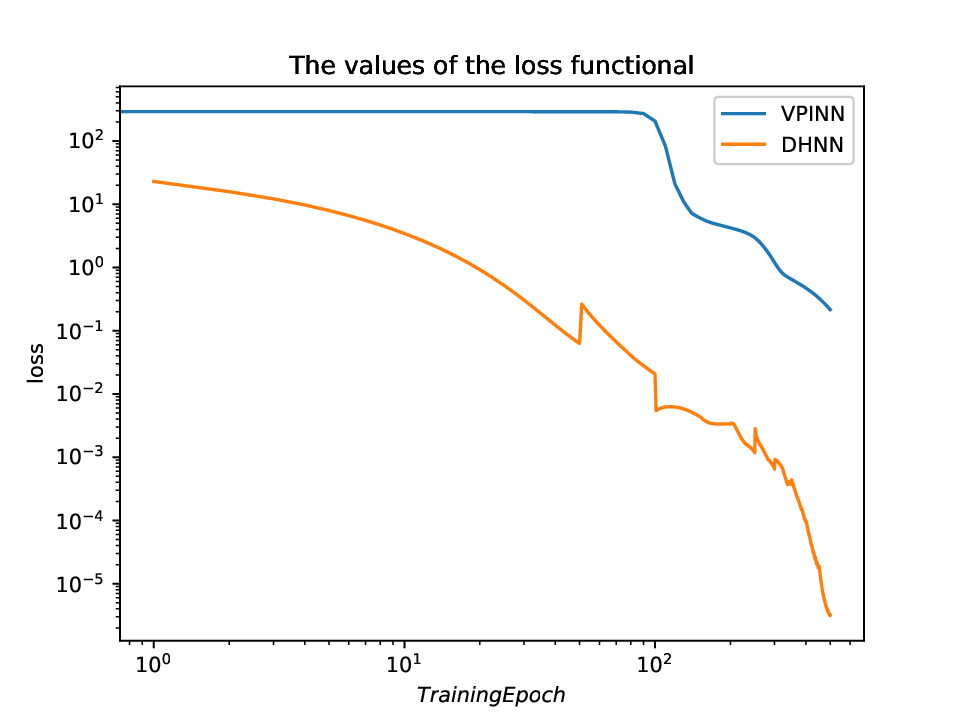}
\end{tabular}
\end{center}
 \caption{The values of the loss functional.}
\label{sin_Possion loss}
\end{figure}

\begin{figure}[H]
\begin{center}
\begin{tabular}{cc}
 \epsfxsize=0.5\textwidth\epsffile{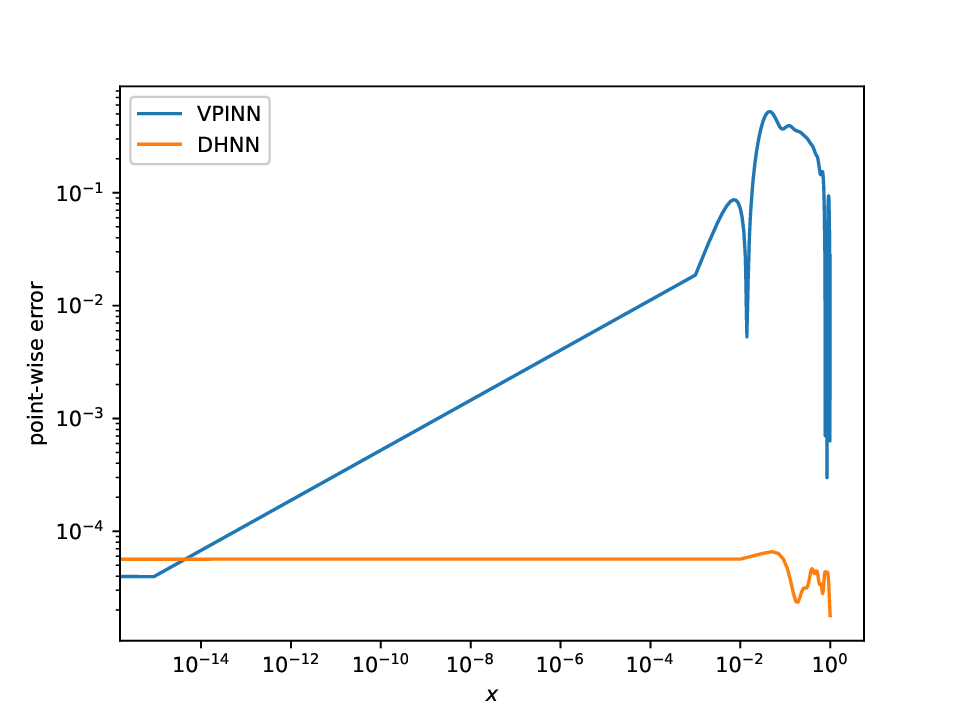}
\end{tabular}
\end{center}
 \caption{The point-wise error.}
\label{sin_Possion error}
\end{figure}

\begin{figure}[H]
  \centering
  \includegraphics[width=0.5\textwidth]{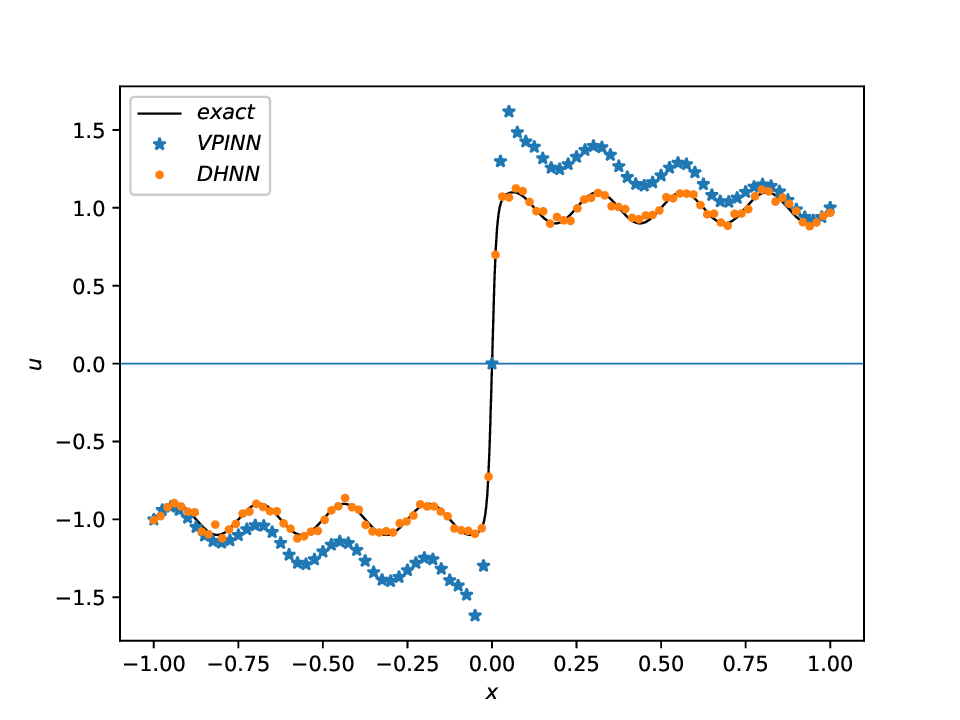} 
  \caption{Comparison between analytical and numerical solutions.}
\label{sin_Possion exact}
\end{figure}

Consider inhomogeneous Possion equation, i.e. $f\neq0$. Let the exact solution be $u = 0.1\sin(4\pi x) + \tanh(5x)$. For DHNN, all nonlinear parameters are initialized to zero. Figures \ref{sin_Possion loss} - \ref{sin_Possion exact} present the following results for both VPINN and DHNN in solving this Poisson equation: the optimization process of loss functionals, point-wise error distributions, and comparisons between exact and numerical solutions. All hidden layers  use $\sin(x)$ as their activation function. The parameters are set as follows:
\begin{align*}
\text{DHNN}: &\quad \ h=\frac{1}{N}=\frac{1}{4},n=10, \alpha=0.04;  \\
\text{VPINN}:    &\quad \ \text{four hidden layers (40 neurons per layer)},\alpha=0.05.
\end{align*}
The experimental results (Figures \ref{sin_Possion loss}-\ref{sin_Possion exact}) demonstrate that the  performance of DHNN is superior over that of VPINN in solving the Poisson equation. Specifically, DHNN converges faster, attaining a loss of $10^{-4}$ within 400 iterations, while the loss of VPINN remains stagnant around 1. Moreover, DHNN exhibits smoother optimization trajectories and yields solutions with more uniform error distributions. Its maximum errors are significantly lower than VPINN's, 
 especially near boundaries and singularities. The numerical solutions also align closely with analytical results and adhere more accurately to boundary conditions. These improvements arise from the neural network architecture of DHNN and advanced optimization algorithm, which together enhance computational efficiency and numerical precision.

\subsection{Helmholtz equations}

Consider Helmholtz equations \cite{Huyuan,yuan} which is formalized, normalizing the wave's velocity to 1, by
\begin{equation}
\left\{ \begin{aligned}
     & -\Delta u(x)-\omega^2u(x) = f  \quad \quad\quad\text{in} \quad \Omega ,\\
    &   (\partial_n+i\omega)u(x) = g   \quad\quad\quad\quad\quad\text{on} \quad\partial\Omega,
                          \end{aligned} \right.
                          \end{equation}
where $g\in H^{\frac{1}{2}}(\partial\Omega)$, $f\in L^2(\Omega)$. The outer normal derivative is referred to by $\partial_n$ and the angular frequency by $\omega$.

\subsubsection{ A homogeneous case}
Consider homogeneous Helmholtz equation, i.e. $f=0$. Let the exact solution be $u=\sin(\omega x)$. The initial values of all nonlinear parameters are set to be 0.

\begin{figure}[H]
\begin{center}
\begin{tabular}{cc}
 \epsfxsize=0.5\textwidth\epsffile{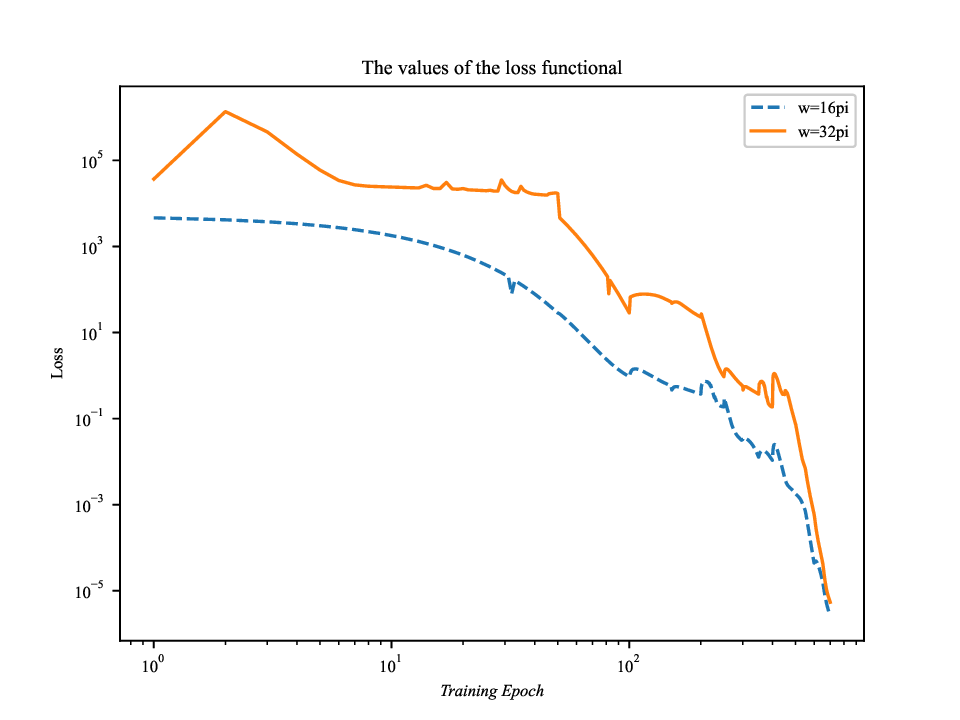}&
 \epsfxsize=0.5\textwidth\epsffile{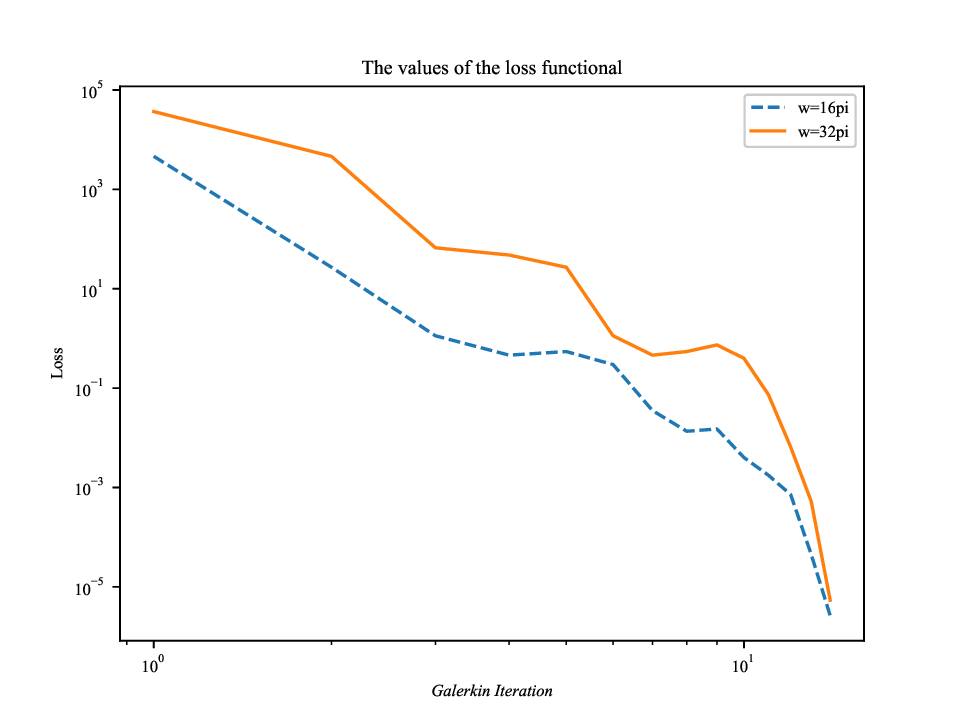}\\
\end{tabular}
\end{center}
 \caption{ The values of the loss functional corresponding to different wave numbers. Left: inner iterations. Right: outer iterations.}
\label{Hel_DWPNN_sig}
\end{figure}

Figure \ref{Hel_DWPNN_sig} shows the optimization process of the loss functional for both inner and outer iterations, with the activation function of the hidden layer being the sigmoid function. The parameters are set  as follows:
\begin{eqnarray*}
\omega=16\pi,h=\frac{1}{N}=\frac{1}{40},n=14, \alpha=0.04;  \\
\omega=32\pi,h=\frac{1}{N}=\frac{1}{60},n=18, \alpha=0.02.
\end{eqnarray*}

From Figure \ref{Hel_DWPNN_sig}, it can be seen that the loss function exhibits a positive convergence trend when the wave numbers $\omega$ are $16\pi$ and $32\pi$, respectively.
In particular, only 14 iterations of the proposed algorithm guarantee the convergence. The value of the loss function decreases by four orders of magnitude in the first three outer iterations and by three orders of magnitude in the last three outer iterations, with a smoother decrease in the middle iterations.

\begin{figure}[H]
\begin{center}
\begin{tabular}{cc}
 \epsfxsize=0.5\textwidth\epsffile{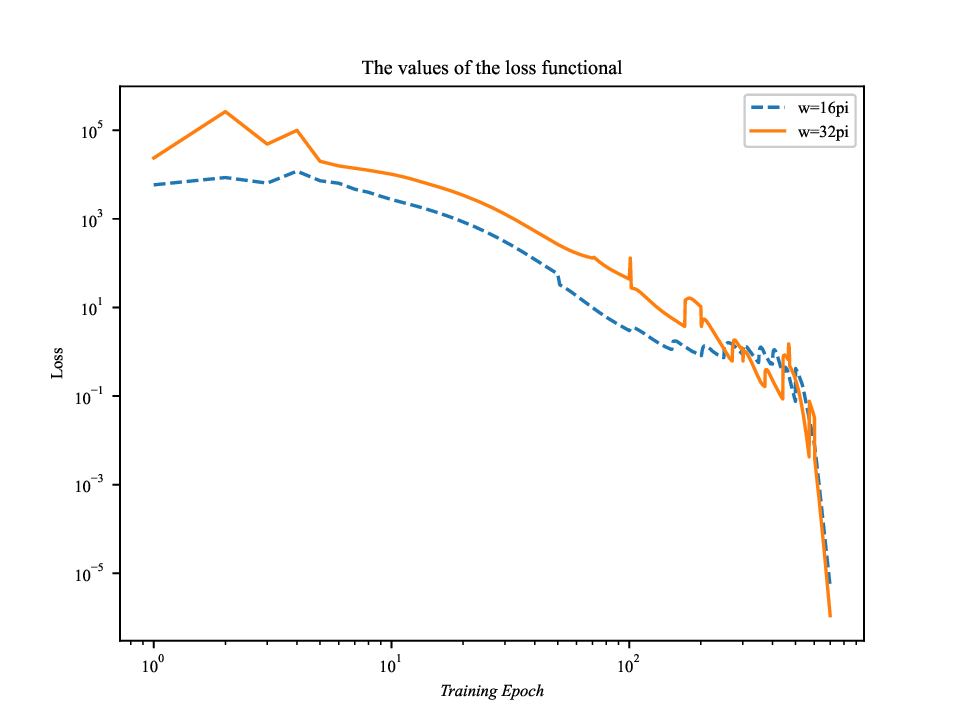}&
 \epsfxsize=0.5\textwidth\epsffile{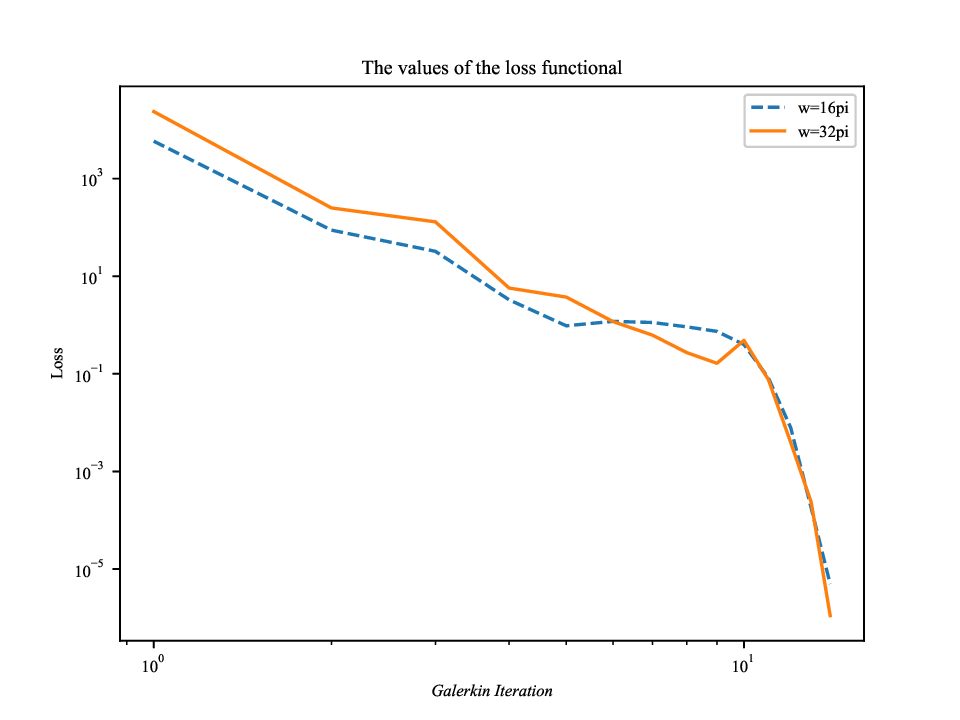}\\
\end{tabular}
\end{center}
 \caption{ The values of the loss functional corresponding to different wave numbers. Left: inner iterations. Right: outer iterations.}
\label{Hel_DWPNN_tanh}
\end{figure}

Figure \ref{Hel_DWPNN_tanh} shows the optimization process of the loss functional for both inner and outer iterations, with the activation function of the hidden layer being the tanh function. The parameters are set as follows:
\begin{eqnarray*}
\omega=16\pi,h=\frac{1}{N}=\frac{1}{40},n=16, \alpha=0.06;  \\
\omega=32\pi,h=\frac{1}{N}=\frac{1}{60},n=20, \alpha=0.04.
\end{eqnarray*}

Figure \ref{Hel_DWPNN_tanh} shows that, although the loss functional experiences occasional slight fluctuations during the iteration process, it still demonstrates an overall convergence trend towards the optimal solution. The value of the loss functional decreased by four orders of magnitude in the first four outer iterations and by two orders of magnitude in the last two outer iterations, with relatively stable decreases in the intermediate outer iterations. This also indicates that the convergence speed of the loss function with a sigmoid activation function is faster compared to the hyperbolic tangent function.

\begin{figure}[H]
\begin{center}
\begin{tabular}{cc}
 \epsfxsize=0.5\textwidth\epsffile{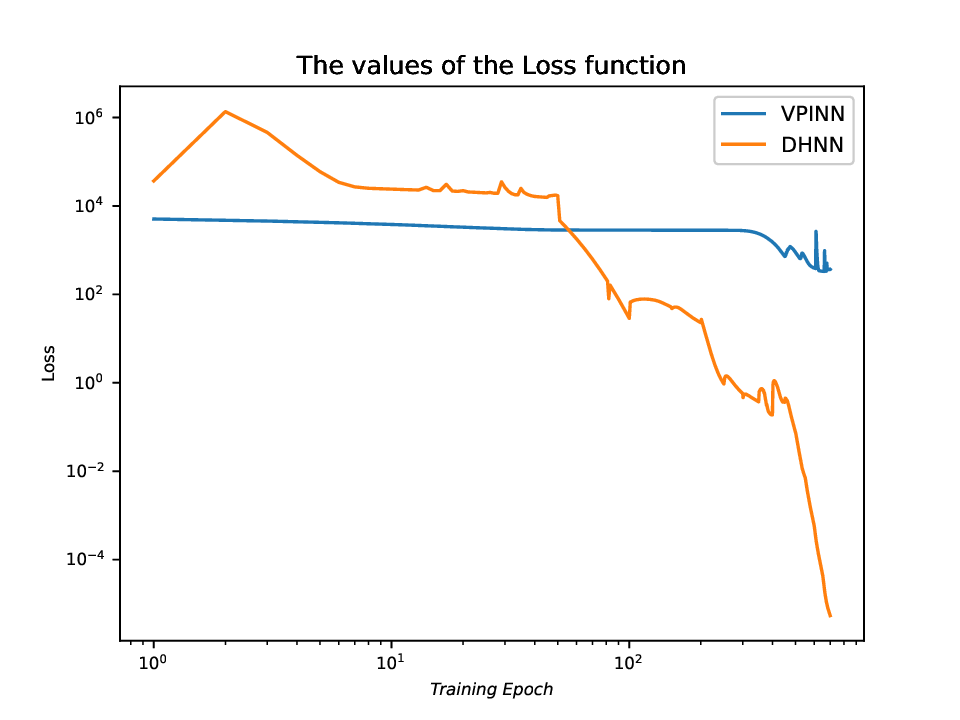}&
 \epsfxsize=0.5\textwidth\epsffile{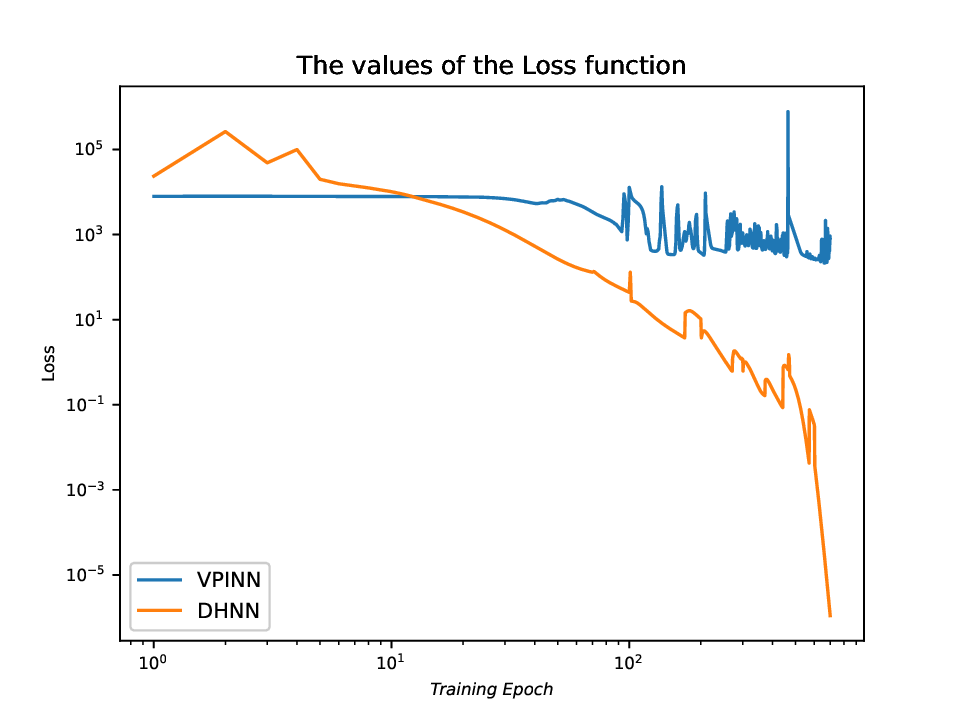}\\
\end{tabular}
\end{center}
 \caption{ Comparison of the loss functional values between the VPINN and the DHNN methods. Left: sigmoid activation function. Right: tanh activation function.}
\label{Hel_VPINN}
\end{figure}

Figure \ref{Hel_VPINN} shows the loss optimization for both VPINN and DHNN methods with sigmoid and tanh activation functions. The parameters of the VPINN method are set as follows:
\begin{align*}
\text{sigmoid}: &\quad \omega=32\pi,\ \text{eight hidden layers (15 neurons per layer)},\ \alpha=0.001;  \\
\text{tanh}:    &\quad \omega=32\pi,\ \text{ten hidden layers (20 neurons per layer)},\ \alpha=0.001.
\end{align*}
During VPINN training, the loss decreases stably but slowly in the early stages, regardless of whether the activation function is sigmoid or tanh. However, after 100 iterations, the loss curve begins oscillating sharply, leading to non-convergence. It is evident that the VPINN cannot converge to the specified accuracy with a limited iteration count.

\subsubsection{An inhomogeneous case}
Taking the exact solution $u=\sin(\sqrt{\omega}x)$, we have $f=(\omega-\omega^2)\sin(\sqrt{\omega}x)$, which indicates the case of $f\neq0$. The initial values for the nonlinear parameters are set as follows: the first $W$ in each element is taken as the midpoint of the current element, while all others $W$ are set to be 0; all parameters $b$ are also set to be 0.

\begin{figure}[H]
\begin{center}
\begin{tabular}{cc}
 \epsfxsize=0.5\textwidth\epsffile{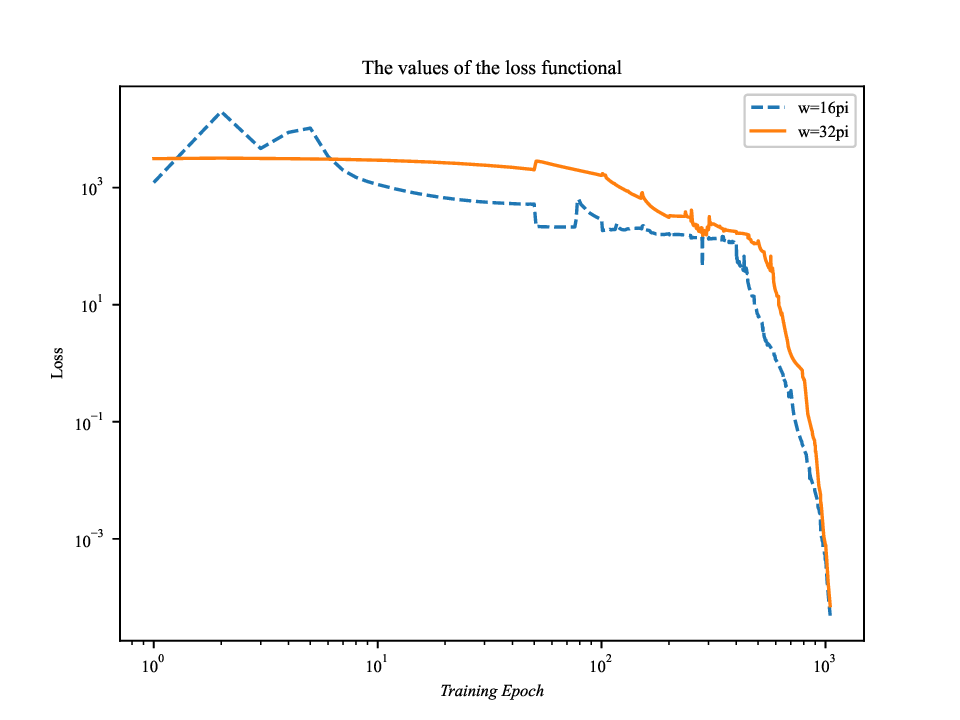}&
 \epsfxsize=0.5\textwidth\epsffile{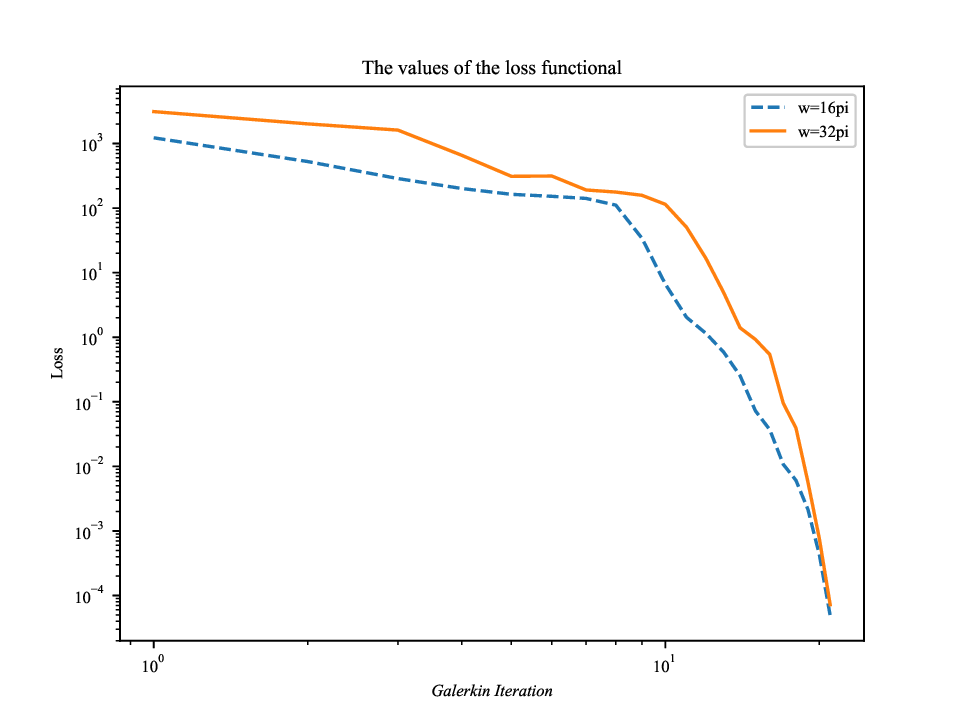}\\
\end{tabular}
\end{center}
 \caption{The values of the loss functional corresponding to different wave numbers. Left: inner iterations. Right: outer iterations.}
\label{inHel_DHNN_sig}
\end{figure}

Figure \ref{inHel_DHNN_sig} shows the optimization process of the loss functional for both inner and outer iterations, with the activation function of the hidden layer being the sigmoid function. The parameters are set as follows:
\begin{eqnarray*}
\omega=16\pi,h=\frac{1}{N}=\frac{1}{30},n=14, \alpha=0.008;  \\
\omega=32\pi,h=\frac{1}{N}=\frac{1}{50},n=16, \alpha=0.006.
\end{eqnarray*}

From Figure \ref{inHel_DHNN_sig}, it can be observed that the number of iterations for the inhomogeneous model is significantly higher compared to the homogeneous model. During the outer iteration process, the decrease in the loss functional value is mainly concentrated in the last few steps, while the initial decline is relatively slow. The outer iteration slightly increases the numerical values of inner iterations when updating linear parameters, although it does not completely enhance them. For example, after 50 inner iterations, the loss functional value may decrease from 0.5 to 0.1; however, an outer iteration step might subsequently increase it back to 0.2, necessitating additional inner iterations to further reduce the loss functional value, which in turn affects the overall iteration speed.

\begin{figure}[H]
\begin{center}
\begin{tabular}{cc}
 \epsfxsize=0.5\textwidth\epsffile{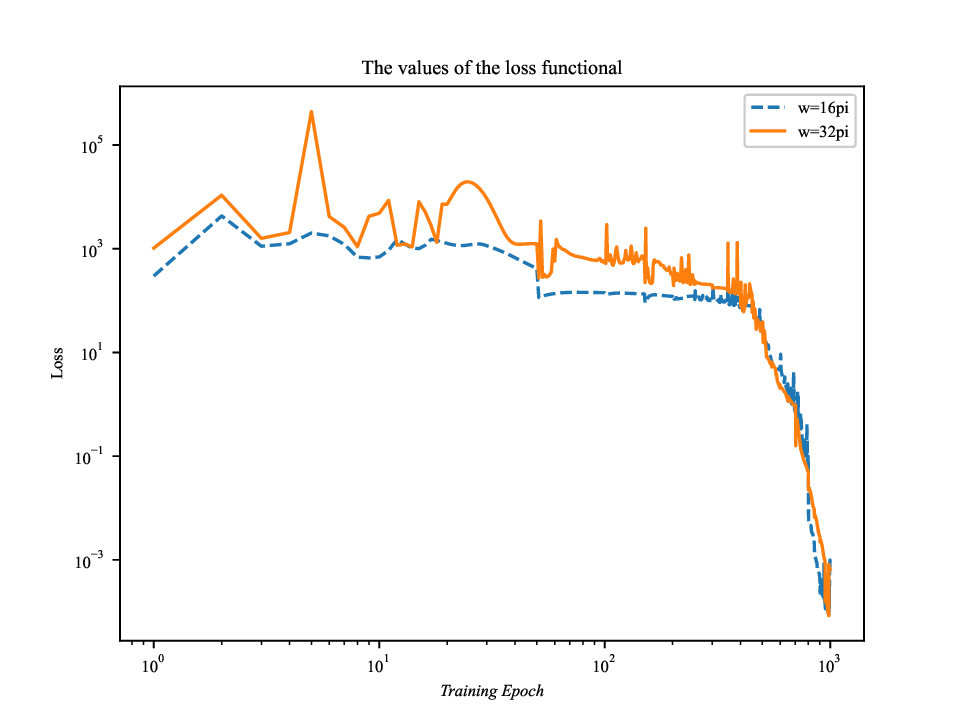}&
 \epsfxsize=0.5\textwidth\epsffile{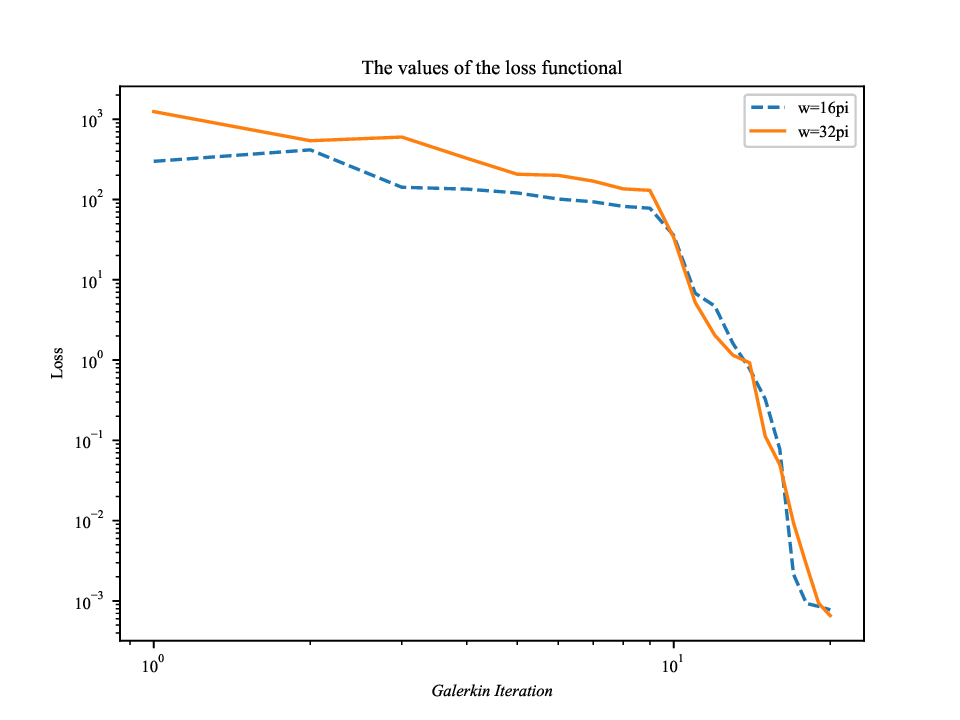}\\
\end{tabular}
\end{center}
 \caption{ The values of the loss functional corresponding to different wave numbers. Left: inner iterations. Right: outer iterations.}
\label{inHel_DHNN_tanh}
\end{figure}

Figure \ref{inHel_DHNN_tanh} shows the optimization process of the loss functional for both inner and outer iterations, with the activation function of the hidden layer being the tanh function. The parameters are set as follows:
\begin{eqnarray*}
\omega=16\pi,h=\frac{1}{N}=\frac{1}{30},n=12, \alpha=0.02; \\
\omega=32\pi,h=\frac{1}{N}=\frac{1}{45},n=14, \alpha=0.01.
\end{eqnarray*}

From Figure \ref{inHel_DHNN_tanh}, it can be observed that during the training process, the loss functional value experiences multiple fluctuations, leading to instability in the early stages of training. However, these short-term fluctuations do not hinder the loss functional value  decreasing to lower levels, indicating that the optimization process is progressing in the right direction. The decrease in the loss functional value during the outer iteration is still primarily concentrated in the last few steps. When using the sigmoid activation function, the algorithm converges to $10^{-5}$; when using the hyperbolic tangent activation function, the algorithm converges to $10^{-4}$. This demonstrates a difference in convergence values between the two activation functions, with the sigmoid activation function converging faster than the tanh function.

\begin{figure}[H]
\begin{center}
\begin{tabular}{cc}
 \epsfxsize=0.5\textwidth\epsffile{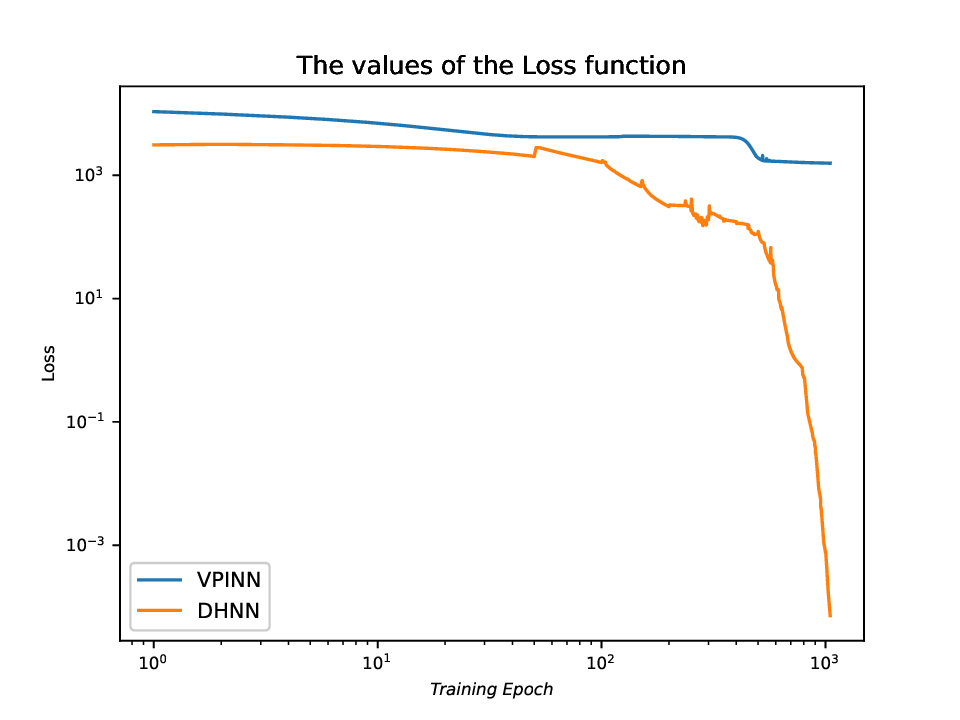}&
 \epsfxsize=0.5\textwidth\epsffile{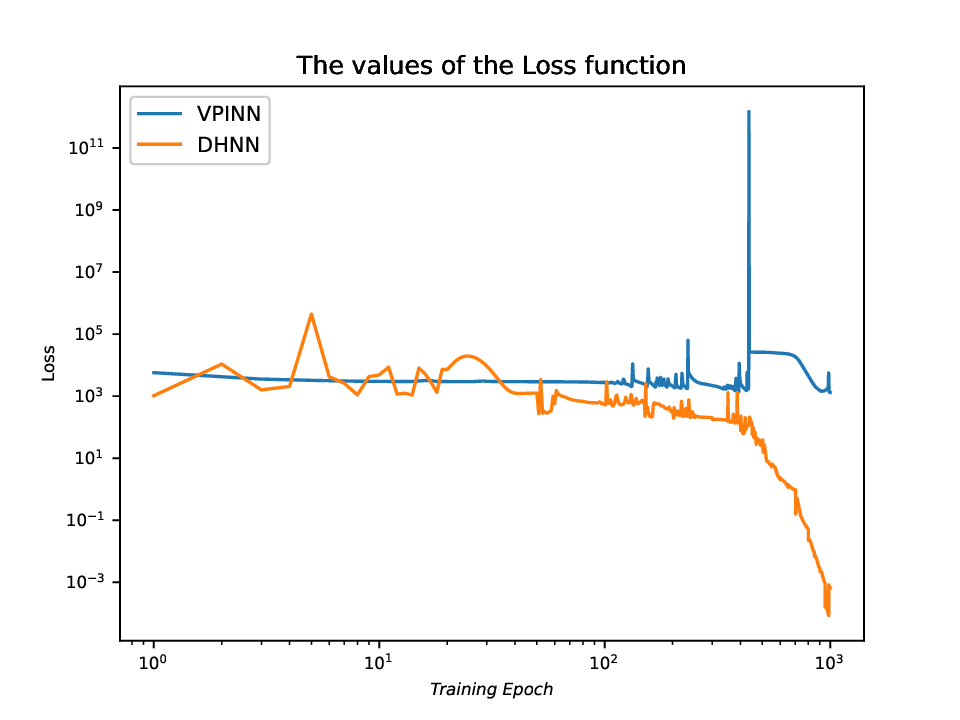}\\
\end{tabular}
\end{center}
 \caption{Comparison of the loss functional values between the VPINN and the DHNN methods. Left: sigmoid activation function. Right: tanh activation function.}
\label{inHel_VPINN}
\end{figure}

Figure \ref{inHel_VPINN} shows the loss optimization for both VPINN and DHNN methods with sigmoid and tanh activation functions. The parameters of the VPINN method are set as follows:
\begin{align*}
\text{sigmoid}: &\quad \omega=32\pi,\ \text{eleven hidden layers (15 neurons per layer)},\ \alpha=0.001;  \\
\text{tanh}:    &\quad \omega=32\pi,\ \text{twelve hidden layers (20 neurons per layer)},\ \alpha=0.001.
\end{align*}
As shown in Figure \ref{inHel_VPINN}, under nonhomogeneous conditions, the VPINN optimization process still exhibits significant numerical instability, characterized by abrupt oscillations in the loss function values that escalate to abnormally high orders of magnitude. Thus, the proposed algorithm demonstrates significantly superior accuracy compared to conventional methods.

\section{Summary}
The discontinuous approximation solutions generated by the training of the discontinuous hybrid neural network proposed in this paper have high accuracy while requiring less computational effort, saving parameter training time. This approach has significant advantages for solving Helmholtz equations with larger wave number. In particular, it adapts well to different equations, thereby opening up an effective new avenue for solving partial differential equations.

\end{document}